# Orthogonal Basis and Motion in Finsler Geometry

Aleks Kleyn

Abstract. Finsler space is differentiable manifold for which Minkowski space is the fiber of the tangent bundle. To understand structure of the reference frame in Finsler space, we need to understand the structure of orthonormal basis in Minkowski space.

In this paper, I considered the definition of orthonormal basis in Minkowski space, the structure of metric tensor relative to orthonormal basis, procedure of orthogonalization. Linear transformation of Minkowski space which preserves the scalar product is called a motion. Linear transformation which maps the orthonormal basis into an infinitely close orthonormal basis is infinitesimal motion. An infinitesimal motion maps orthonormal basis into orthonormal basis.

The set of infinitesimal motions generates Lie algebra, which acts single transitive on basis manifold of Minkowski space. Element of twin representation is called quasimotion of Minkowski space. Quasimotion of event space is called Lorentz transformation.

Contents

## 1. Preface

Geometry is an important tool in modern physics. In particular, research in physics is associated with the search for new geometric constructions. In the paper [2], I showed that for the measurement of spatial and temporal parameters (speed, time delay), the observer uses an orthonormal basis whose vectors are the gauges for measuring time intervals and distances in space. The set of such bases forms a basis manifold. The set of transformations of the basis manifold (Lorentz transformations) generates group, which acts single transitive on the basis manifold.

General relativity does not describe all the phenomena. The goal of string theory and loop gravity is to understand interaction of quantum mechanics and general

Aleks_Kleyn@MailAPS.org.
http://sites.google.com/site/AleksKleyn/.
http://arxiv.org/a/kleyn_a_1.
http://AleksKleyn.blogspot.com/.

relativity. New geometries appeared in frame of these theories. Although these geometries are different from familiar geometry, the set of automorphisms of such geometries generates certain universal algebra $\mathcal{A}$. Our goal is to find a basis of this representation (the definition [3]-3.2.1).

Foundation of the invariance principle is the statement that twin representation of universal algebra $\mathcal{A}$ generates the set of geometrical objects of geometry under consideration.[1] Geometrical object corresponds to the measured physical quantity. The invariance principle is one of the fundamental principles of physics; we use this principle at least since the days of Galileo and Newton. The invariance principle warrants that the measurement of physical quantity with respect to selected basis allows us to predict the measurement of this quantity with respect to another basis. If somebody conducted experiment here and now, then the invariance principle warrants also that we can repeat this experiment in another place and another time.

It was the reason why I was alarmed by statement in paper [5], that violation of Lorentz invariance is possible. I have not found any explanation of this statement. I admit that the form of the Lorentz transformation may be different. I also admit that the transformation of the basis changes and we call it not Lorentz transformation, but differently.

The structure of Finsler geometry is close to the structure of Riemannian geometry. So I attempted to explore the structure of the orthonormal basis and the Lorentz transformations in Finsler geometry. At this stage, since I was interested only in local construction, I explored the orthonormal basis in the Minkowski space.

However, a metric tensor depends on direction. This imposes certain restrictions on the construction and allows us to analyze only infinitely small Lorentz transformation. Infinitesimal Lorentz transformation generates Lie algebra; this gives hope to consider corresponding Lie group. At the same time, even if we find some deviations from the usual structure of Lorentz transformations (for instance, linear transformation maps orthonormal basis into orthonormal one, but it does not preserve the scalar product), this paper gives us opportunity to evaluate the nature of the deviation.

## 2. Finsler Space

In this section, I made definitions similar to definitions in [4].

**Definition 2.1.** A vector space $V$ is called **Minkowski space**[2] when, for the vector space $V$, we define a norm $F$ such that

(1) The norm $F$ is not necessarily positive definite[3]
(2) Function $F(\overline{x})$ is homogeneous of degree 1

$$F(a\overline{x}) = aF(\overline{x}) \quad a > 0 \tag{2.1}$$

[1] I have considered necessary definitions and constructions in sections [3]-3.3, [3]-3.4.

[2] I considered the definition of Minkowski space according to the definition in [6].Although this term calls some association with special relativity, usually it is clear from the context which geometry is referred.

[3] This requirement is due to the fact that we consider applications in general relativity.

(3) Let $\overline{\overline{e}}$ be the basis of vector space $A$. Coordinates of **metric tensor**

$$g_{ij}(\overline{v}) = \frac{\partial^2 F^2(\overline{v})}{\partial v^i \partial v^j} \tag{2.2}$$

form a nonsingular matrix.

□

**Theorem 2.2** (Euler theorem)**.** *Function $f(\overline{x})$, homogeneous of degree $k$,*

$$f(a\overline{x}) = a^k f(\overline{x}) \tag{2.3}$$

*satisfies the differential equation*

$$\frac{\partial f(\overline{x})}{\partial x^i} x^i = k f(\overline{x}) \tag{2.4}$$

*Proof.* We differentiate the equation (2.4) with respect to $a$

$$\frac{df(a\overline{x})}{da} = \frac{da^k}{da} f(\overline{x}) \tag{2.5}$$

According to chain rule, we get

$$\frac{df(a\overline{x})}{da} = \frac{\partial f(a\overline{x})}{\partial a x^i} \frac{dax^i}{da} = \frac{\partial f(a\overline{x})}{\partial a x^i} x^i \tag{2.6}$$

From equations (2.5), (2.6), it follows that

$$\frac{\partial f(a\overline{x})}{\partial a x^i} x^i = k a^{k-1} f(\overline{x}) \tag{2.7}$$

Equation (2.4) follows from equations (2.7) if we assume $a = 1$. □

**Theorem 2.3.** *If $f(\overline{x})$ is homogeneous function of degree $k$, $k > 0$, then partial derivatives $\frac{\partial f(\overline{x})}{\partial x^i}$ are homogeneous functions of degree $k-1$.*

*Proof.* Consider mapping

$$F(x) = \frac{\partial f(\overline{x})}{\partial x^i} x^i \tag{2.8}$$

From equations (2.4), (2.8), it follows that

$$F(\overline{x}) = k f(\overline{x}) \tag{2.9}$$

So, $F(x)$ is ,  $k$. From the equation (2.3), it follows that

$$F(a\overline{x}) = a^k F(\overline{x}) \tag{2.10}$$

From equations (2.8), (2.10), it follows that

$$\frac{\partial f(a\overline{x})}{\partial x^i} a x^i = a^k \frac{\partial f(\overline{x})}{\partial x^i} x^i \tag{2.11}$$

If $k > 0$, then

$$\frac{\partial f(a\overline{x})}{\partial x^i} = a^{k-1} \frac{\partial f(\overline{x})}{\partial x^i}$$

follows from equations (2.11). Therefore, derivatives $\frac{\partial f(\overline{x})}{\partial x^i}$ are homogeneous functions of degree $k-1$. □

**Theorem 2.4.** *The norm of Minkowski space satisfies the differential equations*

$$\frac{\partial F(\overline{a})}{\partial a^i}a^i=F(\overline{a}) \tag{2.12}$$

$$\frac{\partial^2 F(\overline{a})}{\partial a^i\partial a^j}a^i=0 \tag{2.13}$$

$$\frac{1}{2}\frac{\partial^2 F^2(\overline{a})}{\partial a^i\partial a^j}a^ia^j=F^2(\overline{a}) \tag{2.14}$$

*Proof.* The equation (2.12) follows from the statement (2) of the definition 2.1 and the theorem 2.2. According to the theorem 2.2, derivative $\dfrac{\partial F(\overline{x})}{\partial x^i}$ is homogeneous function of degree 0, whence equation (2.13) follows.

Successively differentiating function $F^2$, we get[4]

$$\frac{\partial F^2(\overline{x})}{\partial x^i}=2F(\overline{x})\frac{\partial F(\overline{x})}{\partial x^i}$$

$$\frac{1}{2}\frac{\partial F^2(\overline{x})}{\partial x^j\partial x^i}=\frac{\partial F(\overline{x})}{\partial x^j}\frac{\partial F(\overline{x})}{\partial x^i}+F(\overline{x})\frac{\partial^2 F(\overline{x})}{\partial x^j\partial x^i} \tag{2.15}$$

From equations (2.12), (2.13), (2.15) it follows that

$$\frac{1}{2}\frac{\partial F^2(\overline{x})}{\partial x^j\partial x^i}x^i=\frac{\partial F(\overline{x})}{\partial x^j}\frac{\partial F(\overline{x})}{\partial x^i}x^i+F(\overline{x})\frac{\partial^2 F(\overline{x})}{\partial x^j\partial x^i}x^i=\frac{\partial F(\overline{x})}{\partial x^j}F(\overline{x}) \tag{2.16}$$

From the equation (2.16) it follows that

$$\frac{1}{2}\frac{\partial F^2(\overline{x})}{\partial x^j\partial x^i}x^ix^j=\frac{\partial F(\overline{x})}{\partial x^j}x^jF(\overline{x}) \tag{2.17}$$

The equation (2.14) follows from equations (2.12), (2.17). □

**Theorem 2.5.**

$$\frac{1}{2}g_{ij}(\overline{v})v^iv^j=F^2(\overline{v}) \tag{2.18}$$

*Proof.* The equation (2.18) follows from equations (2.2), (2.14). □

**Theorem 2.6.** *Metric tensor $g_{ij}(\overline{a})$ is homogenious function of degree 0 and satisfies the equation*

$$\frac{\partial g_{ij}(\overline{a})}{\partial a^k}a^k=0 \tag{2.19}$$

*Proof.* From the statement (2) of the definition 2.1, it follows that the mapping $F^2(\overline{v})$ is homogeneous of degree 2. From the theorem 2.3, it follws that function

$$\frac{\partial F^2(\overline{x})}{\partial x^i}$$

is homogeneous of degree 1. From the theorem 2.3 and the definition (2.2), it follows that the function $g_{ij}(\overline{x})$ is homogeneous of degree 0. The equation (2.19) follows from the theorem 2.2. □

**Definition 2.7.** The manifold $M$ is called **Finsler space**, if its tangent space is Minkovsky space and norm $F(x,\overline{v})$ depends continuously on point of tangency $x\in M$. □

[4] see also [6]

*Remark* 2.8. Due to the fact that the norm in the tangent space depends continuously on point of tangency, there is the ability to determine the differential of length of the curve on the manifold

$$dl = F(x, d\overline{x})$$

We usually define Finsler space, and after this we consider tangent to it Minkowski space. In fact, the order of definitions is insignificant. In this paper, Minkowski space is my main object of research. □

## 3. Orthogonality

As Rund noted in [6], there are various definitions of trigonometric functions in Minkowski space. We are primarily interested in the concept of orthogonality.

**Definition 3.1.** Vector $v_1$ **is orthogonal** to vector $v_2$, if

$$g_{ij}(v_1)v_1^i v_2^j = 0$$

□

**Definition 3.2.** The set of vectors $\overline{e}_1$, ..., $\overline{e}_p$ is called orthogonal if

$$\begin{aligned} g_{ij}(\overline{e}_k)e_k^i e_k^j &\neq 0 \\ g_{ij}(\overline{e}_k)e_k^i e_l^j &= 0 \quad k < l \end{aligned} \tag{3.1}$$

The basis $\overline{\overline{e}}$ is called **orthogonal**, if its vectors form orthogonal set. □

**Definition 3.3.** The basis $\overline{\overline{e}}$ is called **orthonormal**, if this is orthogonal basis and its vectors have unit length. □

As we can see from the definition 3.2, relation of orthogonality is noncommutative. Therefore, order of vectors is important when determining an orthogonal basis. There exist different procedures of orthogonalization in Minkowski space. See, for instance, [7], p. 39. Below we consider the orthogonalization procedure proposed in [1], p. 213 - 214.

Any orthogonalization procedure requires positive definite metric. However, if metric is not positive definite, then we can represent Minkowski space as sum of orthogonal spaces $A$ and $B$ such that metric is positive definite in the space $A$, and metric is negative definite in the space $B$. So we also consider orthogonalization procedure in Minkowski space with positive definite metric.

**Theorem 3.4.** *Let* $\overline{e}_1$, ..., $\overline{e}_p$ *be orthogonal set of vectors. Then vectors* $\overline{e}_1$, ..., $\overline{e}_p$ *are linear independent.*

*Proof.* Consider equation

$$a_1\overline{e}_1 + ... + a_p\overline{e}_p = 0 \tag{3.2}$$

From the equation (3.2), it follows that

$$a_1 g_{ij}(\overline{e}_1)e_1^i e_1^j + ... + a_p g_{ij}(\overline{e}_1)e_1^i e_p^j = 0 \tag{3.3}$$

From the condition (3.1) and the equation (3.3) it follows that $a_1 = 0$.

Since we have proved that $a^1 = ... = a^{m-1} = 0$, then the equation (3.2) gets form

$$a_m\overline{e}_m + ... + a_p\overline{e}_p = 0 \tag{3.4}$$

From the equation (3.4), it follows that

$$a_m g_{ij}(\overline{e}_m)e^i_m e^j_m + ... + a_p g_{ij}(\overline{e}_m)e^i_m e^j_p = 0 \tag{3.5}$$

From the condition (3.1) and the equation (3.5) it follows that $a_m = 0$. □

**Theorem 3.5.** *There exists orthonormal basis in Minkowski space.*

*Proof.* Let $n$ be dimension of Minkowski space. Let $\overline{\overline{e}}'$ be a basis in Minkowski space.

We set

$$\overline{e}_1 = \overline{e}'_1$$

Suppose we have defined the set of vectors $\overline{e}_1$, ..., $\overline{e}_m$. In addition we assume that for every $i$, $1 \le i \le m$, the vector $\overline{e}_i$ is linear combination of vectors $\overline{e}'_1$, ..., $\overline{e}'_m$. This assumption also holds for the vector $\overline{e}_{m+1}$, if we represent this vector as

$$\overline{e}_{m+1} = a_1\overline{e}_1 + ... + a_m\overline{e}_m + \overline{e}'_{m+1}$$

$\overline{e}_{m+1} \ne \overline{0}$ because $\overline{\overline{e}}'$ is basis and vector $\overline{e}'_{m+1}$ is not included in the expansion of vectors $\overline{e}_1$, ..., $\overline{e}_m$. For the choice of the vector $\overline{e}'_m$ we require that the vectors $\overline{e}_1$, ..., $\overline{e}_m$. are orthogonal to the vector $\overline{e}_{m+1}$

$$\begin{array}{c} g_{ij}(\overline{e}_1)e^i_1 e^j_{m+1} = 0 \\ ... \\ g_{ij}(\overline{e}_m)e^i_m e^j_{m+1} = 0 \end{array} \tag{3.6}$$

The system of linear equations (3.6) has form

$$\begin{array}{ll} a_1 g_{ij}(\overline{e}_1)e^i_1 e^j_1 & = -g_{ij}(\overline{e}_1)e^i_1 e'^j_{m+1} \\ a_1 g_{ij}(\overline{e}_2)e^i_2 e^j_1 + a_2 g_{ij}(\overline{e}_2)e^i_2 e^j_2 & = -g_{ij}(\overline{e}_2)e^i_2 e'^j_{m+1} \\ ... & \\ a_1 g_{ij}(\overline{e}_m)e^i_m e^j_1 + a_2 g_{ij}(\overline{e}_m)e^i_m e^j_2 + ... + a_m g_{ij}(\overline{e}_m)e^i_m e^j_m & = -g_{ij}(\overline{e}_m)e^i_m e'^j_{m+1} \end{array}$$

Therefore, the solution of the system of linear equations (3.6) has form

$$a_1 = -\frac{g_{ij}(\overline{e}_1)e^i_1 e'^j_{m+1}}{g_{ij}(\overline{e}_1)e^i_1 e^j_1}$$

$$a_2 = -\frac{g_{ij}(\overline{e}_2)e^i_2 e'^j_{m+1} + a_1 g_{ij}(\overline{e}_2)e^i_2 e^j_1}{g_{ij}(\overline{e}_2)e^i_2 e^j_2}$$

...

$$a_m = -\frac{g_{ij}(\overline{e}_m)e^i_m e'^j_{m+1} - \sum_{k=1}^{m-1} a_k g_{ij}(\overline{e}_m)e^i_m e^j_k}{g_{ij}(\overline{e}_m)e^i_m e^j_m}$$

By continuing this process, we obtain an orthogonal basis $\overline{\overline{e}}$. We can normalize vectors of the basis $\overline{\overline{e}}$ according to rule

$$E_k = g_{ij}(\overline{e}_k)e^i_k e^j_k$$

$$\overline{e}_k \rightarrow \frac{1}{\sqrt{E_k}}\overline{e}_k$$

$$k = 1, ..., n$$

□

**Theorem 3.6.** *Let $\overline{\overline{e}}$ be orthonormal basis of Minkowski space. If we write the coordinates of the metric tensor $g_{ij}(\overline{e}_k)$ relative to basis $\overline{\overline{e}}$ as matrix*

$$\begin{pmatrix} g_{11}(\overline{e}_1) & ... & g_{1n}(\overline{e}_1) \\ ... & ... & ... \\ g_{n1}(\overline{e}_n) & ... & g_{nn}(\overline{e}_n) \end{pmatrix} \tag{3.7}$$

*then the matrix (3.7) is a triangular matrix whose diagonal elements are* 1*.*

*Proof.* Evidently, we can assume $e^i_j = \delta^i_j$. According to the definition 3.2, 3.3

$$\begin{aligned} g_{ij}(\overline{e}_k)\delta^i_k\delta^j_k &= 1 \\ g_{ij}(\overline{e}_k)\delta^i_k\delta^j_l &= 0 \quad k < l \end{aligned} \tag{3.8}$$

From equations (3.8), it follows that

$$\begin{aligned} g_{kk}(\overline{e}_k) &= 1 \\ g_{kl}(\overline{e}_k) &= 0 \quad k < l \end{aligned}$$

Therefore, $g_{kl}(\overline{e}_k)$ is arbitrary, when $k > l$. □

## 4. Motion of Minkowski space

The structure of Minkowski space is close to the structure Euclidean space. Automorphism of Minkowski space as well as automorphism of Euclidean space is called **motion**.

Since the orthogonality relation is not symmetric, then the structure of metric tensor in an orthonormal basis changes; this creates a dramatic expansion of the set of orthonormal bases. In particular, since scalar product

$$g_{ij}(\overline{e}_k)e^i_k e^j_l \quad k > l$$

is not defined, then we cannot require that automorphism of Minkowski space preserves scalar product.

According to the theorem [3]-3.2.6, the motion maps orthonormal basis $\overline{\overline{e}}$ into orthonormal basis $\overline{\overline{e}}'$

$$\begin{aligned} g_{ij}(\overline{e}_k)e^i_k e^j_k &= 1 \quad & g_{ij}(\overline{e}'_k)e'^i_k e'^j_k &= 1 \\ g_{ij}(\overline{e}_k)e^i_k e^j_l &= 0 \quad & g_{ij}(\overline{e}'_k)e'^i_k e'^j_l &= 0 \quad k < l \end{aligned} \tag{4.1}$$

Since the motion does not change the basis, then mappings $g_{ij}$ also do not change. However, the argument of the mapping $g_{ij}$ changes. Therefore, we can consider only an infinitesimal motion.

**Theorem 4.1.** *Let infinitesimal motion*

$$a'^i = a^j(\delta^i_j + A^i_j dt)$$

*map orthonormal basis $\overline{\overline{e}}$ to orthonormal basis $\overline{\overline{e}}'$*

$$e'^i_k = e^j_k(\delta^i_j + A^i_j dt) \tag{4.2}$$

*Then* $(k \le l)$

$$g_{kp}(\overline{e}_k)A^p_l+g_{pl}(\overline{e}_k)A^p_k+\left.\frac{\partial g_{kl}(\overline{a})}{\partial a^p}\right|_{\overline{a}=\overline{e}_k}A^p_k=0 \tag{4.3}$$

*Proof.* Coordinates of metric tensor change according to rule

$$g_{ij}(\overline{v}')=g_{ij}(\overline{v})+\frac{\partial g_{ij}(\overline{v})}{\partial v^k}A^k_l v^l dt \tag{4.4}$$

For $k\le l$, from equations (4.2), (4.4), it follows that $(\overline{v}=\overline{e}'_k)$

$$\begin{aligned}
&g_{ij}(\overline{e}'_k){e'}^i_k{e'}^j_l\\
&=(g_{ij}(\overline{e}_k)+\left.\frac{\partial g_{ij}(\overline{a})}{\partial a^p}\right|_{\overline{a}=\overline{e}_k}e^r_kA^p_rdt)(e^i_k+e^m_kA^i_mdt)(e^j_l+e^p_lA^j_pdt)\\
&=(g_{ij}(\overline{e}_k)+\left.\frac{\partial g_{ij}(\overline{a})}{\partial a^p}\right|_{\overline{a}=\overline{e}_k}e^r_kA^p_rdt)(e^i_ke^j_l+e^i_ke^j_pA^p_ldt+e^i_mA^m_kdte^j_l)\\
&=g_{ij}(\overline{e}_k)e^i_ke^j_l+g_{ij}(\overline{e}_k)e^i_ke^j_pA^p_ldt+g_{ij}(\overline{e}_k)e^i_mA^m_kdte^j_l\\
&+\left.\frac{\partial g_{ij}(\overline{a})}{\partial a^p}\right|_{\overline{a}=\overline{e}_k}e^r_kA^p_rdte^i_ke^j_l
\end{aligned} \tag{4.5}$$

From equations (4.1), (4.5), it follows that $(k\le l)$

$$0=g_{ij}(\overline{e}_k)e^i_ke^j_pA^p_ldt+g_{ij}(\overline{e}_k)e^i_mA^m_kdte^j_l+\left.\frac{\partial g_{ij}(\overline{a})}{\partial a^p}\right|_{\overline{a}=\overline{e}_k}e^r_kA^p_rdte^i_ke^j_l \tag{4.6}$$

Since $e^i_k=\delta^i_k$, then the equation (4.3) follows from equation (4.6). □

**Theorem 4.2.** *The product of infinitesimal motions of Mincovsky space is infinitesimal motion of Mincovsky space.*

*Proof.* Let

$$f(\overline{a})^i=a^k(\delta^i_k+A^i_kdt)$$
$$g(\overline{a})^i=a^k(\delta^i_k+B^i_kdt)$$

be infinitesimal motions of Mincovsky space. Transformation $fg$ has form

$$\begin{aligned}
f(g(\overline{a}))^i&=f(a^l(\delta^i_l+B^i_ldt))\\
&=a^l(\delta^i_l+B^i_ldt)(\delta^l_k+A^l_kdt)\\
&=a^l(\delta^i_k+A^i_kdt+B^i_kdt)\\
&=a^l(\delta^i_k+(A^i_k+B^i_k)dt)
\end{aligned}$$

From the theorem 4.1, it follows that coordinates of mappings $f$ and $g$ satisfy equation

$$\begin{aligned}
&g_{kp}(\overline{e}_k)A^p_l+g_{pl}(\overline{e}_k)A^p_k+\left.\frac{\partial g_{kl}(\overline{a})}{\partial a^p}\right|_{\overline{a}=\overline{e}_k}A^p_k=0\\
&g_{kp}(\overline{e}_k)B^p_l+g_{pl}(\overline{e}_k)B^p_k+\left.\frac{\partial g_{kl}(\overline{a})}{\partial a^p}\right|_{\overline{a}=\overline{e}_k}B^p_k=0
\end{aligned}\qquad k\le l \tag{4.7}$$

From the equation (4.7), it follows that $(k\le l)$

$$
\begin{aligned}
&g_{kp}(\overline{e}_k)(A_l^p+B_l^p)+g_{pl}(\overline{e}_k)(A_k^p+B_k^p)+\left.\frac{\partial g_{kl}(\overline{a})}{\partial a^p}\right|_{\overline{a}=\overline{e}_k}(A_k^p+B_k^p)\\
=&g_{kp}(\overline{e}_k)A_l^p+g_{pl}(\overline{e}_k)A_k^p+\left.\frac{\partial g_{kl}(\overline{a})}{\partial a^p}\right|_{\overline{a}=\overline{e}_k}A_k^p\\
+&g_{kp}(\overline{e}_k)B_l^p+g_{pl}(\overline{e}_k)B_k^p+\left.\frac{\partial g_{kl}(\overline{a})}{\partial a^p}\right|_{\overline{a}=\overline{e}_k}B_k^p\\
=&0
\end{aligned}
$$

Therefore, the mapping $fg$ is infinitesimal motion of Mincovsky space. □

## 5. Quasimotion of Minkovsky Space

Linear transformation of the basis

$$\overline{e}'_i=A_i^j\overline{e}_j$$

is called **quasimotion of Minkovsky space**.

Coordinates of vector

$$\overline{a}=a^i\overline{e}_i$$

transform according to the rule

$$a'^j=A^{-1}{}^j_i a^i \tag{5.1}$$

From equations (2.18), (5.1) it follows that

$$g_{ij}(\overline{a})a^ia^j=g'_{kl}(\overline{a})a'^ka'^l=g'_{kl}(\overline{a})A^{-1}{}^k_i a^iA^{-1}{}^l_j a^j \tag{5.2}$$

From the equation (5.2), it follows that

$$g'_{kl}(\overline{a})=g_{ij}(\overline{a})A_k^iA_l^j \tag{5.3}$$

Therefore, $g_{ij}(\overline{a})$ is tensor.

Consider infinitesimal transformation

$$\overline{a}'=\overline{a}+d\overline{a}$$

Then there is infinitesimal transformation of metric tensor

$$g_{ij}(\overline{a}')=g_{ij}(\overline{a})+\frac{\partial g_{ij}(\overline{a})}{\partial a^k}da^k \tag{5.4}$$

Consider infinitesimal quasimotion

$$\overline{e}'_i=\overline{e}_j(\delta_i^j+A_i^jdt)$$

According to (5.3), (5.4), it follows that there is infinitesimal transformation of metric tensor[5]

$$
\begin{aligned}
g'_{kl}(\overline{e}'_p)&=g_{ij}(\overline{e}'_p)(\delta_k^i\delta_l^j+\delta_k^iA_l^jdt+A_k^i\delta_l^jdt)\\
&=(g_{ij}(\overline{e}_p)+\left.\frac{\partial g_{ij}(\overline{a})}{\partial a^m}\right|_{\overline{a}=\overline{e}_p}e_r^mA_p^rdt)(\delta_k^i\delta_l^j+\delta_k^iA_l^jdt+A_k^i\delta_l^jdt)\\
&=(g_{ij}(\overline{e}_p)+\left.\frac{\partial g_{ij}(\overline{a})}{\partial a^m}\right|_{\overline{a}=\overline{e}_p}A_p^mdt)(\delta_k^i\delta_l^j+\delta_k^iA_l^jdt+A_k^i\delta_l^jdt)\\
&=g_{kl}(\overline{e}_p)+g_{ij}(\overline{e}_p)(\delta_k^iA_l^jdt+A_k^i\delta_l^jdt)+\left.\frac{\partial g_{kl}(\overline{a})}{\partial a^m}\right|_{\overline{a}=\overline{e}_p}A_p^mdt
\end{aligned}
\tag{5.5}
$$

[5] We use equation $e_i^m=\delta_i^m$.

For $k\le l$,

$$g_{kl}(\overline{e}_k)=g'_{kl}(\overline{e}'_k) \tag{5.6}$$

From equations (5.5), (5.6), it follows that $(k\le l)$

$$g_{kj}(\overline{e}_k)A^j_l+g_{il}(\overline{e}_k)A^i_k+\left.\frac{\partial g_{kl}(\overline{a})}{\partial a^m}\right|_{\overline{a}=\overline{e}_k}A^m_p=0$$

Therefore, the set of infinitesimal quasimotions of Minkovsky space generates the same algebra as the set of infinitesimal motions of Minkovsky space.

# 7. Index

# Ортогональный базис и движение в финслеровой геометрии

Александр Клейн

Аннотация. Финслерово пространство - это дифференцируемое многообразие, для которого пространство Минковского является слоем касательного расслоения. Для того, чтобы понять строение системы отсчёта в Финслеровом пространстве, мы должны понять структуру ортонормального базиса в пространстве Минковского.

В статье рассмотрено определение ортонормального базиса в пространстве Минковского, структура метрического тензора относительно ортонормального базиса, процедура ортогонализации. Линейное преобразование пространства Минковского, сохраняющее скалярное произведение, называется движением. Линейное преобразование отображающее ортонормальный базис в бесконечно близкий ортонормальный базис, является инфинитезимальным движением. Инфинитезимальное движение отображает ортонормальный базис в ортонормальный базис.

Множество инфинитезимальных движений порождает алгебру Ли, однотранзитивно действующую на многообразии базисов пространства Минковского. Элемент парного представления называется квазидвижением пространства Минковского. Квазидвижение пространства событий, называется преобразованием Лоренца.

## Содержание

## 1. Предисловие

Геометрия является важным инструментом в современной физике. В частности, исследование в физике связано с поиском новых геометрических конструкций. В статье [2], я показал, что при измерении пространственно временных параметров (скорость, задержка времени) наблюдатель пользуется ортонормальным базисом, векторы которого являются калибрами для измерения интервалов времени и расстояний в пространстве. Множество таких базисов

Aleks_Kleyn@MailAPS.org.
http://sites.google.com/site/AleksKleyn/.
http://arxiv.org/a/kleyn_a_1.
http://AleksKleyn.blogspot.com/.

формирует многообразие базисов. Множество преобразований многообразия базисов (преобразований Лоренца) порождает группу, действующую однотранзитивно на многообразии базисов.

Общая теория относительности не описывает все явления. Задача теории струн и петлевой гравитации - понять как взаимодействуют квантовая механика и общая теория относительности. В рамках этих теорий появились новые геометрии. Хотя эти геометрии отличаются от привычной нам геометрии, множество автоморфизмов этих геометрий порождает некоторую универсальную алгебру $\mathcal{A}$. Наша задача найти базис этого представления (определение [3]-3.2.1).

Основа принципа инвариантности состоит в том, что парное представление универсальной алгебры $\mathcal{A}$ порождает множество геометрических объектов рассматриваемой геометрии.[1] Геометрический объект соответствует измеряемой физической величине. Принцип инвариантности является одним из фундаментальных принципов физики, и мы пользуемся этим принципом, по крайней мере начиная со времён Галилея и Ньютона. Принцип инвариантности гарантирует, что измерение физической величины относительно данного базиса позволяет предсказать измерение этой величины относительно другого базиса. Принцип инвариантности гарантирует также, что эксперимент, проведенный здесь и сейчас, можно повторить в другом месте и в другое время.

Поэтому меня насторожило утверждение в статье [5], что возможно нарушение инвариантности Лоренца. Никаких разъяснений этого утверждения я не нашёл. Я допускаю, что форма преобразований Лоренца может быть другой. Я также допускаю, что закон преобразования базиса изменится настолько, что мы будем называть его не преобразованием Лоренца, а как-то иначе.

Структура финслеровой геометрии близка структуре римановой геометрии. Поэтому я предпринял попытку исследовать структуру ортонормального базиса и преобразования Лоренца в финслеровой геометрии. Так как на этом этапе меня интересовали только локальные построения, то я ограничился изучением ортонормального базиса в пространстве Минковского.

Однако метрический тензор зависит от направления. Это накладывает определённые ограничения на построения и позволяет проанализировать только бесконечно малые преобразования Лоренца. Бесконечно малые преобразования Лоренца порождают алгебру Ли, и это даёт надежду рассмотреть соответствующую группу Ли. В тоже время, даже если мы обнаружим некоторые отклонения от привычной структуры преобразований Лоренца (например, линейное преобразование отображает ортонормальный базис в ортонормальный, но не сохраняет скалярное произведение), то эта статья даёт возможность оценить характер отклонения.

## 2. Финслерово пространство

Определения в этом разделе даны по аналогии с определениями в [4].

[1] Я рассмотрел необходимые определения и построения в разделах [3]-3.3, [3]-3.4.

**Определение 2.1.** Векторное пространство $V$ называется **пространством Минковского**,[2] если в векторном пространстве $V$ определена норма $F$ такая, что

(1) Норма $F$ не обязательно положительно определена[3]
(2) Функция $F(\overline{x})$ однородна степени 1

$$F(a\overline{x}) = aF(\overline{x}) \quad a > 0 \tag{2.1}$$

(3) Пусть $\overline{\overline{e}}$ - базис векторного пространства $A$. Координаты **метрического тензора**

$$g_{ij}(\overline{v}) = \frac{\partial^2 F^2(\overline{v})}{\partial v^i \partial v^j} \tag{2.2}$$

порождают невырожденную матрицу.

□

**Теорема 2.2** (Теорема Эйлера)**.** *Функция $f(\overline{x})$, однородная степени $k$,*

$$f(a\overline{x}) = a^k f(\overline{x}) \tag{2.3}$$

*удовлетворяет дифференциальному уравнению*

$$\frac{\partial f(\overline{x})}{\partial x^i} x^i = kf(\overline{x}) \tag{2.4}$$

*Доказательство.* Продифференцируем равенство (2.4) по $a$

$$\frac{df(a\overline{x})}{da} = \frac{da^k}{da} f(\overline{x}) \tag{2.5}$$

Согласно правилу дифференцирования по частям, мы имеем

$$\frac{df(a\overline{x})}{da} = \frac{\partial f(a\overline{x})}{\partial ax^i} \frac{dax^i}{da} = \frac{\partial f(a\overline{x})}{\partial ax^i} x^i \tag{2.6}$$

Из равенств (2.5), (2.6) следует

$$\frac{\partial f(a\overline{x})}{\partial ax^i} x^i = ka^{k-1} f(\overline{x}) \tag{2.7}$$

Равенство (2.4) следует из равенства (2.7) если положить $a = 1$. □

**Теорема 2.3.** *Если $f(\overline{x})$ - функция, однородная степени $k$, $k > 0$, то частные производные $\dfrac{\partial f(\overline{x})}{\partial x^i}$ являются функциями, однородными степени $k-1$.*

*Доказательство.* Рассмотрим отображение

$$F(x) = \frac{\partial f(\overline{x})}{\partial x^i} x^i \tag{2.8}$$

Из равенств (2.4), (2.8) следует

$$F(\overline{x}) = kf(\overline{x}) \tag{2.9}$$

---

[2] Я рассмотрел определение пространства Минковского согласно определению в [6], с. 28 - 32. Хотя этот термин вызывает некоторые ассоциации со специальной теорией относительности, обычно из контекста ясно о какой геометрии идёт речь.

[3] Это требование связано с тем, что мы рассматриваем приложения в общей теории относительности.

Следовательно, $F(x)$ - функция, однородная степени $k$. Из равенства (2.3) следует

$$F(a\overline{x}) = a^k F(\overline{x}) \tag{2.10}$$

Из равенств (2.8), (2.10) следует

$$\frac{\partial f(a\overline{x})}{\partial x^i} a x^i = a^k \frac{\partial f(\overline{x})}{\partial x^i} x^i \tag{2.11}$$

Если $k > 0$, то из равенства (2.11) следует

$$\frac{\partial f(a\overline{x})}{\partial x^i} = a^{k-1} \frac{\partial f(\overline{x})}{\partial x^i}$$

Следовательно, производные $\dfrac{\partial f(\overline{x})}{\partial x^i}$ являются однородными функциями степени $k-1$. □

**Теорема 2.4.** *Норма пространства Минковского удовлетворяет дифференциальным уравнениям*

$$\frac{\partial F(\overline{a})}{\partial a^i} a^i = F(\overline{a}) \tag{2.12}$$

$$\frac{\partial^2 F(\overline{a})}{\partial a^i \partial a^j} a^i = 0 \tag{2.13}$$

$$\frac{1}{2} \frac{\partial^2 F^2(\overline{a})}{\partial a^i \partial a^j} a^i a^j = F^2(\overline{a}) \tag{2.14}$$

*Доказательство.* Равенство (2.12) следует из утверждения (2) определения 2.1 и теоремы 2.2. Согласно теореме 2.2 производная $\dfrac{\partial F(\overline{x})}{\partial x^i}$ является однородной функцией степени 0, откуда следует равенство (2.13).

Последовательно дифференцируя функцию $F^2$, мы получим[4]

$$\frac{\partial F^2(\overline{x})}{\partial x^i} = 2F(\overline{x}) \frac{\partial F(\overline{x})}{\partial x^i}$$

$$\frac{1}{2} \frac{\partial F^2(\overline{x})}{\partial x^j \partial x^i} = \frac{\partial F(\overline{x})}{\partial x^j} \frac{\partial F(\overline{x})}{\partial x^i} + F(\overline{x}) \frac{\partial^2 F(\overline{x})}{\partial x^j \partial x^i} \tag{2.15}$$

Из равенств (2.12), (2.13), (2.15) следует

$$\frac{1}{2} \frac{\partial F^2(\overline{x})}{\partial x^j \partial x^i} x^i = \frac{\partial F(\overline{x})}{\partial x^j} \frac{\partial F(\overline{x})}{\partial x^i} x^i + F(\overline{x}) \frac{\partial^2 F(\overline{x})}{\partial x^j \partial x^i} x^i = \frac{\partial F(\overline{x})}{\partial x^j} F(\overline{x}) \tag{2.16}$$

Из равенства (2.16) следует

$$\frac{1}{2} \frac{\partial F^2(\overline{x})}{\partial x^j \partial x^i} x^i x^j = \frac{\partial F(\overline{x})}{\partial x^j} x^j F(\overline{x}) \tag{2.17}$$

Равенство (2.14) следует из равенств (2.12), (2.17). □

**Теорема 2.5.**

$$\frac{1}{2} g_{ij}(\overline{v}) v^i v^j = F^2(\overline{v}) \tag{2.18}$$

*Доказательство.* Равенство (2.18) является следствием равенств (2.2), (2.14). □

[4] см. также [6], с. 24.

**Теорема 2.6.** *Метрический тензор $g_{ij}(\overline{a})$ является однородной функцией степени 0 и удовлетворяет уравнению*

$$\frac{\partial g_{ij}(\overline{a})}{\partial a^k}a^k = 0 \tag{2.19}$$

*Доказательство.* Из утверждения (2) определения 2.1 следует, что отображение $F^2(\overline{v})$ однородно степени 2. Из теоремы 2.3 следует, что функция

$$\frac{\partial F^2(\overline{x})}{\partial x^i}$$

однородна степени 1. Из теоремы 2.3 и определения (2.2) следует, что функция $g_{ij}(\overline{x})$ однородна степени 0. Равенство (2.19) следует из теоремы 2.2. □

**Определение 2.7.** Многообразие $M$ называется **финслеровым пространством**, если его касательное пространство являются пространством Минковского и норма $F(x,\overline{v})$ непрерывно зависит от точки касания $x \in M$. □

*Замечание* 2.8. Следствием того, что норма в касательном пространстве непрерывно зависит от точки касания, является возможность определения дифференциала длины кривой на многообразии

$$dl = F(x, d\overline{x})$$

Обычно сперва определяют финслерово пространство, а потом рассматривают касательное к нему пространство Минковского. На самом деле порядок определений несущественен. В этой статье, моим основным объектом исследования является пространство Минковского. □

## 3. Ортогональность

Как отметил Рунд в [6], с. 47, существуют различные определения тригонометрических функций в пространстве Минковского. Нас интересует прежде всего понятие ортогональности.

**Определение 3.1.** Вектор $v_1$ **ортогонален** вектору $v_2$, если

$$g_{ij}(v_1)v_1^i v_2^j = 0$$

□

**Определение 3.2.** Множество векторов $\overline{e}_1$, ..., $\overline{e}_p$ называется ортогональным, если

$$\begin{aligned} &g_{ij}(\overline{e}_k)e_k^i e_k^j \neq 0 \\ &g_{ij}(\overline{e}_k)e_k^i e_l^j = 0 \quad k < l \end{aligned} \tag{3.1}$$

Базис $\overline{\overline{e}}$ называется **ортогональным**, если его векторы формируют оргональное множество. □

**Определение 3.3.** Базис $\overline{\overline{e}}$ называется **ортонормированым**, если это ортогональный базис и его векторы имеют единичную длину. □

Как мы видим из определения 3.2, отношение ортогональности некоммутативно. Следовательно, порядок векторов важен при определении ортогонального базиса. Существуют различные процедуры ортогонализации в пространстве Минковского. Смотри, например, [7], с. 39. Ниже мы рассмотрим процедуру ортогонализации, предложенную в [1], с. 213 - 214.

Любая процедура ортогонализации предполагает положительно определённую метрику. Однако, если метрика не является положительно определённой, то пространство Минковского можно представить в виде суммы ортогональных пространств $A$ и $B$ таких, что в пространстве $A$ метрика положительно определена, а в пространстве $B$ метрика отрицательно определена. Поэтому мы также будем рассматривать процедуру ортогонализации в пространстве Минковского с положительно определённой метрикой.

**Теорема 3.4.** *Пусть* $\overline{e}_1, ..., \overline{e}_p$ *- ортогональное множество векторов. Тогда векторы* $\overline{e}_1, ..., \overline{e}_p$ *линейно независимы.*

*Доказательство.* Рассмотрим равенство

$$a_1\overline{e}_1 + ... + a_p\overline{e}_p = 0 \tag{3.2}$$

Из равенства (3.2) следует

$$a_1 g_{ij}(\overline{e}_1)e_1^i e_1^j + ... + a_p g_{ij}(\overline{e}_1)e_1^i e_p^j = 0 \tag{3.3}$$

Из условия (3.1) и равенства (3.3) следует $a_1 = 0$.

Если мы доказали, что $a^1 = ... = a^{m-1} = 0$, то равенство (3.2) примет вид

$$a_m\overline{e}_m + ... + a_p\overline{e}_p = 0 \tag{3.4}$$

Из равенства (3.4) следует

$$a_m g_{ij}(\overline{e}_m)e_m^i e_m^j + ... + a_p g_{ij}(\overline{e}_m)e_m^i e_p^j = 0 \tag{3.5}$$

Из условия (3.1) и равенства (3.5) следует $a_m = 0$. □

**Теорема 3.5.** *Ортонормированый базис в пространстве Минковского существует.*

*Доказательство.* Пусть $n$ - размерность пространства Минковского. Пусть $\overline{\overline{e}}'$ - базис в пространстве Минковского.

Мы положим

$$\overline{e}_1 = \overline{e}'_1$$

Допустим мы построили множество векторов $\overline{e}_1, ..., \overline{e}_m$. Дополнительно предположим, что для всякого $i$, $1 \le i \le m$, вектор $\overline{e}_i$ является линейной комбинацией векторов $\overline{e}'_1, ..., \overline{e}'_m$. Это предположение будет выполнено и для вектора $\overline{e}_{m+1}$, если мы этот вектор представим в виде

$$\overline{e}_{m+1} = a_1\overline{e}_1 + ... + a_m\overline{e}_m + \overline{e}'_{m+1}$$

$\overline{e}_{m+1} \ne \overline{0}$ так как $\overline{\overline{e}}'$ - базис, и вектор $\overline{e}'_{m+1}$ не входит в разложение векторов $\overline{e}_1, ..., \overline{e}_m$. Для выбора вектора $\overline{e}'_m$ мы потребуем, чтобы векторы $\overline{e}_1, ..., \overline{e}_m$. были ортогональны вектору $\overline{e}_{m+1}$

$$\begin{array}{c} g_{ij}(\overline{e}_1)e_1^i e_{m+1}^j = 0 \\ ... \\ g_{ij}(\overline{e}_m)e_m^i e_{m+1}^j = 0 \end{array} \tag{3.6}$$

Система линейных уравнений (3.6) имеет вид

$$\begin{aligned}
&a_1 g_{ij}(\overline{e}_1)e_1^i e_1^j && = -g_{ij}(\overline{e}_1)e_1^i e'^j_{m+1}\\
&a_1 g_{ij}(\overline{e}_2)e_2^i e_1^j + a_2 g_{ij}(\overline{e}_2)e_2^i e_2^j && = -g_{ij}(\overline{e}_2)e_2^i e'^j_{m+1}\\
&...\\
&a_1 g_{ij}(\overline{e}_m)e_m^i e_1^j + a_2 g_{ij}(\overline{e}_m)e_m^i e_2^j + ... + a_m g_{ij}(\overline{e}_m)e_m^i e_m^j && = -g_{ij}(\overline{e}_m)e_m^i e'^j_{m+1}
\end{aligned}$$

Следовательно, решение системы линейных уравнений (3.6) имеет вид

$$\begin{aligned}
a_1 &= -\frac{g_{ij}(\overline{e}_1)e_1^i e'^j_{m+1}}{g_{ij}(\overline{e}_1)e_1^i e_1^j}\\
a_2 &= -\frac{g_{ij}(\overline{e}_2)e_2^i e'^j_{m+1} + a_1 g_{ij}(\overline{e}_2)e_2^i e_1^j}{g_{ij}(\overline{e}_2)e_2^i e_2^j}\\
&...\\
a_m &= -\frac{g_{ij}(\overline{e}_m)e_m^i e'^j_{m+1} - \sum_{k=1}^{m-1} a_k g_{ij}(\overline{e}_m)e_m^i e_k^j}{g_{ij}(\overline{e}_m)e_m^i e_m^j}
\end{aligned}$$

Продолжая этот процесс, мы получим ортогональный базис $\overline{\overline{e}}$. Мы можем нормировать векторы базиса $\overline{\overline{e}}$ согласно правилу

$$\begin{aligned}
E_k &= g_{ij}(\overline{e}_k)e_k^i e_k^j\\
\overline{e}_k &\rightarrow \frac{1}{\sqrt{E_k}}\overline{e}_k\\
k &= 1, ..., n
\end{aligned}$$

□

**Теорема 3.6.** *Пусть $\overline{\overline{e}}$ - ортонормальный базис пространства Минковского. Если мы запишем координаты метрического тензора $g_{ij}(\overline{e}_k)$ относительно базиса $\overline{\overline{e}}$ в виде матрицы*

$$\begin{pmatrix}
g_{11}(\overline{e}_1) & ... & g_{1n}(\overline{e}_1)\\
... & ... & ...\\
g_{n1}(\overline{e}_n) & ... & g_{nn}(\overline{e}_n)
\end{pmatrix} \tag{3.7}$$

*то матрица (3.7) является треугольной матрицей, диагональные элементы которой равны* 1.

*Доказательство.* Очевидно, мы можем положить $e_j^i = \delta_j^i$. Согласно определениям 3.2, 3.3

$$\begin{aligned}
&g_{ij}(\overline{e}_k)\delta_k^i \delta_k^j = 1\\
&g_{ij}(\overline{e}_k)\delta_k^i \delta_l^j = 0 \quad k < l
\end{aligned} \tag{3.8}$$

Из равенств (3.8) следует

$$\begin{aligned}
&g_{kk}(\overline{e}_k) = 1\\
&g_{kl}(\overline{e}_k) = 0 \quad k < l
\end{aligned}$$

Следовательно, $g_{kl}(\overline{e}_k)$ произвольно, если $k>l$. □

## 4. Движение пространства Минковского

Структура пространства Минковского близка структуре евклидова пространства. Автоморфизм пространства Минковского, так же как и автоморфизм евклидова пространства, называется **движением**.

Так как отношение ортогональности не симметрично, то это приводит к изменению структуры метрического тензора в ортонормальном базисе и резкому расширению множества ортонормальных базисов. В частности, так как скалярное произведение

$$g_{ij}(\overline{e}_k)e^i_k e^j_l \quad k>l$$

не определено, то мы не можем требовать, что автоморфизм пространства Минковского сохраняет скалярное произведение.

Согласно теореме [3]-3.2.6, движение отображает ортонормальный базис $\overline{\overline{e}}$ в ортонормальный базис $\overline{\overline{e}}'$

$$\begin{aligned} &g_{ij}(\overline{e}_k)e^i_k e^j_k=1 \quad g_{ij}(\overline{e}'_k)e'^i_k e'^j_k=1\\ &g_{ij}(\overline{e}_k)e^i_k e^j_l=0 \quad g_{ij}(\overline{e}'_k)e'^i_k e'^j_l=0 \quad k<l \end{aligned} \tag{4.1}$$

Поскольку при движении базис пространства не меняется, то отображения $g_{ij}$ также не меняются. Однако аргумент отображения $g_{ij}$ меняется. Поэтому мы можем рассмотреть только инфинитезимальное движение.

**Теорема 4.1.** *Пусть бесконечно малое движение*

$$a'^i=a^j(\delta^i_j+A^i_j dt)$$

*отображает ортонормированный базис* $\overline{\overline{e}}$ *в ортонормированный базис* $\overline{\overline{e}}'$

$$e'^i_k=e^j_k(\delta^i_j+A^i_j dt) \tag{4.2}$$

*Тогда* $(k\le l)$

$$g_{kp}(\overline{e}_k)A^p_l+g_{pl}(\overline{e}_k)A^p_k+\left.\frac{\partial g_{kl}(\overline{a})}{\partial a^p}\right|_{\overline{a}=\overline{e}_k}A^p_k=0 \tag{4.3}$$

*Доказательство.* Координаты метрического тензора меняются согласно правилу

$$g_{ij}(\overline{v}')=g_{ij}(\overline{v})+\frac{\partial g_{ij}(\overline{v})}{\partial v^k}A^k_l v^l dt \tag{4.4}$$

Для $k\le l$ из равенств (4.2), (4.4), следует $(\overline{v}=\overline{e}'_k)$

$$\begin{aligned} &g_{ij}(\overline{e}'_k)e'^i_k e'^j_l\\ &=(g_{ij}(\overline{e}_k)+\left.\frac{\partial g_{ij}(\overline{a})}{\partial a^p}\right|_{\overline{a}=\overline{e}_k}e^r_k A^p_r dt)(e^i_k+e^m_k A^i_m dt)(e^j_l+e^p_l A^j_p dt)\\ &=(g_{ij}(\overline{e}_k)+\left.\frac{\partial g_{ij}(\overline{a})}{\partial a^p}\right|_{\overline{a}=\overline{e}_k}e^r_k A^p_r dt)(e^i_k e^j_l+e^i_k e^j_p A^p_l dt+e^i_m A^m_k dt e^j_l)\\ &=g_{ij}(\overline{e}_k)e^i_k e^j_l+g_{ij}(\overline{e}_k)e^i_k e^j_p A^p_l dt+g_{ij}(\overline{e}_k)e^i_m A^m_k dt e^j_l\\ &+\left.\frac{\partial g_{ij}(\overline{a})}{\partial a^p}\right|_{\overline{a}=\overline{e}_k}e^r_k A^p_r dt e^i_k e^j_l \end{aligned} \tag{4.5}$$

Из равенств (4.1), (4.5), следует $(k \le l)$

$$0 = g_{ij}(\overline{e}_k)e^i_k e^j_p A^p_l dt + g_{ij}(\overline{e}_k)e^i_m A^m_k dt e^j_l + \left.\frac{\partial g_{ij}(\overline{a})}{\partial a^p}\right|_{\overline{a}=\overline{e}_k} e^r_k A^p_r dt e^i_k e^j_l \tag{4.6}$$

Так как $e^i_k = \delta^i_k$, то равенство (4.3) следует из равенства (4.6). □

**Теорема 4.2.** *Произведение инфинитезимальных движений пространства Минковского является инфинитезимальным движением пространства Минковского.*

*Доказательство.* Пусть

$$f(\overline{a})^i = a^k(\delta^i_k + A^i_k dt)$$
$$g(\overline{a})^i = a^k(\delta^i_k + B^i_k dt)$$

инфинитезимальные движения пространства Минковского. Преобразование $fg$ имеет вид

$$\begin{aligned} f(g(\overline{a}))^i &= f(a^l(\delta^i_l + B^i_l dt)) \\ &= a^l(\delta^i_l + B^i_l dt)(\delta^l_k + A^l_k dt) \\ &= a^l(\delta^i_k + A^i_k dt + B^i_k dt) \\ &= a^l(\delta^i_k + (A^i_k + B^i_k)dt) \end{aligned}$$

Из теоремы 4.1, следует, что координаты отображений $f$ и $g$ удовлетворяют равенствам

$$\begin{aligned} g_{kp}(\overline{e}_k)A^p_l + g_{pl}(\overline{e}_k)A^p_k + \left.\frac{\partial g_{kl}(\overline{a})}{\partial a^p}\right|_{\overline{a}=\overline{e}_k} A^p_k &= 0 \\ g_{kp}(\overline{e}_k)B^p_l + g_{pl}(\overline{e}_k)B^p_k + \left.\frac{\partial g_{kl}(\overline{a})}{\partial a^p}\right|_{\overline{a}=\overline{e}_k} B^p_k &= 0 \end{aligned} \qquad k \le l \tag{4.7}$$

Из равенства (4.7) следует $(k \le l)$

$$\begin{aligned} &g_{kp}(\overline{e}_k)(A^p_l + B^p_l) + g_{pl}(\overline{e}_k)(A^p_k + B^p_k) + \left.\frac{\partial g_{kl}(\overline{a})}{\partial a^p}\right|_{\overline{a}=\overline{e}_k} (A^p_k + B^p_k) \\ =&g_{kp}(\overline{e}_k)A^p_l + g_{pl}(\overline{e}_k)A^p_k + \left.\frac{\partial g_{kl}(\overline{a})}{\partial a^p}\right|_{\overline{a}=\overline{e}_k} A^p_k \\ +&g_{kp}(\overline{e}_k)B^p_l + g_{pl}(\overline{e}_k)B^p_k + \left.\frac{\partial g_{kl}(\overline{a})}{\partial a^p}\right|_{\overline{a}=\overline{e}_k} B^p_k \\ =&0 \end{aligned}$$

Следовательно, отображение $fg$ является инфинитезимальным движением пространства Минковского. □

## 5. Квазидвижение пространства Минковского

Линейное преобразование многообразия базисов

$$\overline{e}'_i = A^j_i \overline{e}_j$$

называется **квазидвижение пространства Минковского**.

Координаты вектора

$$\overline{a} = a^i \overline{e}_i$$

преобразуются согласно правилу

$$a'^j = A^{-1\cdot}{}^j_i a^i \tag{5.1}$$

Из равенств (2.18), (5.1) следует

$$g_{ij}(\overline{a})a^i a^j = g'_{kl}(\overline{a})a'^k a'^l = g'_{kl}(\overline{a})A^{-1\cdot}{}^k_i a^i A^{-1\cdot}{}^l_j a^j \tag{5.2}$$

Из равенства (5.2) следует

$$g'_{kl}(\overline{a}) = g_{ij}(\overline{a})A^i_k A^j_l \tag{5.3}$$

Следовательно, $g_{ij}(\overline{a})$ является тензором.

Рассмотрим бесконечно малое преобразование

$$\overline{a}' = \overline{a} + d\overline{a}$$

Тогда метрический тензор испытывает бесконечно малое преобразование

$$g_{ij}(\overline{a}') = g_{ij}(\overline{a}) + \frac{\partial g_{ij}(\overline{a})}{\partial a^k} da^k \tag{5.4}$$

Рассмотрим бесконечно малое квазидвижение

$$\overline{e}'_i = \overline{e}_j(\delta^j_i + A^j_i dt)$$

Согласно (5.3), (5.4), следует, что метрический тензор испытывает бесконечно малое преобразование[5]

$$\begin{aligned} g'_{kl}(\overline{e}'_p) &= g_{ij}(\overline{e}'_p)(\delta^i_k\delta^j_l + \delta^i_k A^j_l dt + A^i_k \delta^j_l dt) \\ &= (g_{ij}(\overline{e}_p) + \left.\frac{\partial g_{ij}(\overline{a})}{\partial a^m}\right|_{\overline{a}=\overline{e}_p} e^m_r A^r_p dt)(\delta^i_k\delta^j_l + \delta^i_k A^j_l dt + A^i_k \delta^j_l dt) \\ &= (g_{ij}(\overline{e}_p) + \left.\frac{\partial g_{ij}(\overline{a})}{\partial a^m}\right|_{\overline{a}=\overline{e}_p} A^m_p dt)(\delta^i_k\delta^j_l + \delta^i_k A^j_l dt + A^i_k \delta^j_l dt) \\ &= g_{kl}(\overline{e}_p) + g_{ij}(\overline{e}_p)(\delta^i_k A^j_l dt + A^i_k \delta^j_l dt) + \left.\frac{\partial g_{kl}(\overline{a})}{\partial a^m}\right|_{\overline{a}=\overline{e}_p} A^m_p dt \end{aligned} \tag{5.5}$$

При $k \le l$,

$$g_{kl}(\overline{e}_k) = g'_{kl}(\overline{e}'_k) \tag{5.6}$$

Из равенств (5.5), (5.6), следует, что $(k \le l)$

$$g_{kj}(\overline{e}_k)A^j_l + g_{il}(\overline{e}_k)A^i_k + \left.\frac{\partial g_{kl}(\overline{a})}{\partial a^m}\right|_{\overline{a}=\overline{e}_k} A^m_p = 0$$

Следовательно, множество бесконечно малых квазидвижений пространства Минковского порождает туже алгебру, что и множество бесконечно малых движений пространства Минковского.

[5]Мы используем равенство $e^m_i = \delta^m_i$.

## 6. Список литературы

# 7. Предметный указатель